\documentclass [11pt]{article}

\newtheorem{definition}{Definition}
\newtheorem{proposition}{Proposition}
\newtheorem{condition}{Condition}
\newtheorem{theorem}{Theorem}
\newtheorem{lemma}[theorem]{Lemma}

\newtheorem{remark}{Remark}
\usepackage{amsmath}
\usepackage{amsfonts}
\usepackage{amssymb}
\usepackage{verbatim}

\title{A Study on the Amount of Random Graph Groupies}

\author{
Daodi Lu
\thanks{
Department of Mathematics, California Institute of Technology, Pasadena, CA 91125, USA.
dlu@caltech.edu}
}

\begin{document}

\maketitle

\begin{abstract}
    In 1980, Ajtai, Komlos and Szemer{\'e}di defined ``groupie'':
    Let $G=(V,E)$ be a simple graph, $|V|=n$, $|E|=e$.
    For a vertex $v\in V$, let $r(v)$ denote the sum of the degrees of the vertices adjacent to $v$.
    We say $v\in V$ is a {\it groupie}, if
    $\frac{r(v)}{\deg(v)}\geq\frac{e}{n}.$
    In this paper, we prove that in random graph $B(n,p)$, $0<p<1$, the proportion of groupies
    converges in probability towards $\Phi(1)\approx0.8413$ as $n$ approaches infinity, where
    $\Phi(x)$ is the distribution function of standard normal distribution $N(0,1)$.
    We also discuss the asymptotic behavior of the proportion of groupies in complete bipartite graph $B(n_1,n_2,p)$.
\end{abstract}

\section{Introduction}   The definition of groupie was first given in \cite{AKS1980}.
    In that paper, Ajtai, Komlos and Szemer{\'e}di used the fact that every nonempty simple graph has at least one groupie
    to give an upper bound for Ramsey number $R(3,k)$.

\begin{definition}\label{groupie}
    Suppose $G=(V,E)$ is a nonempty simple graph, with $|V|=n$ and $|E|=e$. For a vertex $v$ of $G$,
    denote by $r(v)$ the sum of the degrees of neighbors of $v$. The vertex $v$ is called a groupie,
    if the average degree of the neighbors of $v$ is not less than the average degree of all vertices in $G$, i.e.
    \begin{equation}
        \frac{r(v)}{\deg(v)}\geq\frac{e}{n}.
    \end{equation}
    For the case that $v$ is isolated, $v$ is a groupie if and only if all vertices in $G$ are isolated.
\end{definition}
    The references \cite{BEHST1994}, \cite{PST1995}, \cite{M1996}, and \cite{H2007} discuss the properties of groupies in simple graphs.
    Mackey \cite{M1996} proved that there are at least two groupies in any simple graphs with at least two vertices.
    The original definition of groupie was revised in \cite{VT2009} and \cite{S2010}. In these papers, the term
    ``groupie'' is defined as a vertex whose degree is not less than the average degree of its neighbors in a simple graph.
    Fernandez de la Vega and Tuza \cite{VT2009} proved that the proportion of groupie in graph $B(n,p)$
    is almost always very near to $\frac{1}{2}$ as $n\rightarrow\infty$.
    Shang \cite{S2010} investigated the amount of groupies in bipartite graph $B(n_1,n_2,p)$, and proved that the proportion of groupie in
    $B(n_1,n_2,p)$ is almost always very near to $\frac{1}{2}$ when $n_1,n_2\rightarrow \infty$ and $n_1=n_2$.
    However, he did not show whether the proportion of groupies converges in probability when $n_1,n_2$ go to infinity together in general case.
    In our paper, the term ``groupie'' is used as Definition \ref{groupie}.
    In Section \ref{B(n,p)}, we will investigate the proportion of groupie defined in Definition \ref{groupie} in a random graph $B(n,p)$.
    As the number of vertices $n\rightarrow\infty$, the proportion of groupie converges in probability towards $\Phi(1)\approx0.8413$,
    where $\Phi(x)$ is the cumulative distribution function of the standard normal distribution $N(0,1)$. The main result in this section is:

\begin{theorem} \label{main}
    Suppose $G=B(n,p)$, $0<p<1$, is a complete random simple graph on $n$ vertices. Let $N(n)$ be the number of groupies in $G$.
    Then for any $\epsilon>0$, we have
    \begin {equation}
        P \bigg( \bigg|\frac{N(n)}{n}-\Phi(1) \bigg|>\epsilon \bigg)\rightarrow0, \quad \textrm{as }n\rightarrow\infty,
    \end{equation}
    where $\Phi(x)$ is the cumulative distribution function of the standard normal distribution $N(0,1)$.
\end{theorem}

    In Section \ref{B(n_1,n_2,p)}, we will discuss the asymptotic behavior of the proportion of groupie in a random bipartite graphs $B(n_1,n_2,p)$
    as $n_1,n_2\rightarrow\infty$ and the ratio $\frac{n_1}{n_2}\rightarrow \alpha$, where $\alpha$ is a fixed nonnegative number.
    If $|n_1-n_2|\rightarrow \infty$ as $n_1,n_2\rightarrow\infty$, the proportion of groupies in $B(n_1,n_2,p)$ converges in probability towards $\frac{\max(1,\alpha)}{1+\alpha}$,
    while when $\alpha=1$ and $n_1-n_2=c$ where $c$ is a fixed integer, the proportion of groupies converges in probability towards $\frac{1}{2}(\Phi(1+\frac{pc}{2(1-p)})+\Phi(1-\frac{pc}{2(1-p)}))$.
    This limit is monotone decreasing with the absolute value of $c$.
    When $c=0$, the limit is $\Phi(1)\approx 0.8413$.
    As $|c|$ is large, it is near to $\frac{1}{2}$,
    which coincides with the limit when $\alpha=1$ and $|n_1-n_2|\rightarrow\infty$ as $n_1,n_2\rightarrow\infty$.

\section{In The Complete Random Graph $B(n,p)$}\label{B(n,p)}

In random graph $G=B(n,p)$, denote the $n$ vertices to be $v_1,v_2,\ldots,v_n$.
For $k=1,2,\ldots,n$, let $A_k$ be the event ``The vertex $v_k$ is a groupie''.
Let $X_k$ be the characteristic function of $A_k$, i.e. $X_k=1_{A_k}$.\\

We claim that we can exclude the case that there exists an isolated vertex in graph $G$, because
as $n\rightarrow\infty$,
\begin{displaymath}
P(\textrm{there exists an isolated vertex in }B(n,p))  \leq n\cdot(1-p)^{n-1} \rightarrow 0.
\end{displaymath}
We first prove:

\begin{lemma}\label{EX_1}
The limit $\lim_{n\rightarrow\infty}E[X_1]-\Phi(1)$ is equal to $0$.
\end{lemma}

Note that $E[X_1]=P(A_1)$, i.e. $E[X_1]$ is the probability that the vertex $v_1$ is a groupie.
Suppose $i$ is the degree of $v_1$.
In a graph $G=B(n,p)$, let $V_1$ be the set of $i$ vertices which are adjacent to $v_1$;
$V_2$ the set of $n-1-i$ vertices which are not adjacent to $v_1$;
$e_1$ the number of edges whose two vertices are all in the vertex set $V_1$;
$e_2$ the number of edges whose two vertices are all in the vertex set $V_2$;
$e_3$ the number of edges with one vertex in the vertex set $V_1$, and the other vertex in $V_2$.
Then the event $A_1$ is equivalent to $$\frac{2(i+e_1+e_2+e_3)}{n}\leq \frac{i+2e_1+e_3}{i},$$
i.e. $$2(n-i)e_1+(n-2i)(e_3+i)-2ie_2\geq0.$$
Note that the conditional expectation
\begin{align*}
&E[2(n-i)e_1+(n-2i)(e_3+i)-2ie_2;i]\\
=&ip((n-i)(i-1)+(n-2i)(n-1-i)-(n-1-i)(n-2-i))+(n-2i)i\\
=&i((n-2)p+(n-2i)),
\end{align*}
and the conditional variance
\begin{align*}
    &Var[2(n-i)e_1+(n-2i)(e_3+i)-2ie_2;i]\\
    =&\bigg((2(n-i))^2\frac{i(i-1)}{2}+(n-2i)^2 i(n-1-i)\\
     &+(2i)^2\frac{(n-1-i)(n-2-i)}{2}\bigg)(p-p^2)\\
    \triangleq & \sigma^2.
\end{align*}

Since $i\thicksim B(n-1,p)$, for any $\epsilon_0>0$,
we can pick $N_0\in \mathbb{N}$,
such that for all $n>N_0$, $P(|i-pn|\geq n^{0.5}\Omega(n))<\epsilon_0$,
where $\Omega(n)$ approaches infinity as $n\rightarrow\infty$, with the speed slower than any power of $n$.
Denote by $F_1(x)$, $F_2(x)$, $F_3(x)$ the cumulative distribution function of the binomial random variables $\frac{2(n-i)e_1}{\sigma}$, $\frac{-2ie_2}{\sigma}$, $\frac{(n-2i)e_3}{\sigma}$,
and let $\Phi_1(x)$, $\Phi_2(x)$, $\Phi_3(x)$ be the cumulative distribution function of the normal distribution in central limit theorem corresponding to these binomial random variables. Suppose $X, X_1,\ldots,X_n$ are independent identically distributed $d$-dimensional random vectors,
and $X=(X^{(1)},\ldots,X^{(d)})$. In addition, suppose the third moment of all components of $X$ exist, and the first moment $E[X]=0$.
Let $$Y=\frac{X_1+\ldots+X_n}{n},$$
$F(x_1,\ldots,x_d)$ the cumulative distribution function of $Y$,
and $\Phi(x_1,\ldots,x_d)$ the cumulative distribution function of the $d$-dimensional normal distribution $\widetilde{Y}$ that has the same first moment and covariance matrix as $Y$.
By the Berry-Esseen Theorem \cite{B1941}\cite{E1942}\cite{B1945}\cite{C1948},
there exist a constant $C(d)$ that only depends on the dimension $d$, so that
$$\sup_{x\in \mathbb{R}^d}|F(x)-\Phi(x)|\leq\frac{C(d)}{n^{\frac{1}{2}}}\sum_{i=1}^{d}\frac{E[|X^{(i)}|^3]}{(E[(X^{(i)})^2])^{\frac{3}{2}}}.$$

By the 1-dimensional Berry-Esseen theorem\cite{B1941}\cite{E1942}, there exists a constant $C$ that does not depend on $n,i,p$, so that
$$\sup_{x\in \mathbb{R}}|F_j(x)-\Phi_j(x)|<\frac{C(p^2+(1-p)^2)}{\sqrt{p(1-p)}}\frac{1}{\sqrt{E_j}},\mspace{15mu} j=1,2,3,$$
where $E_j$ denote the number of random edges which influence $e_j$, after given the degree $i$ of vertex $v_1$, for example, $E_1=\binom{i}{2}$.

If $|i-pn|<n^{0.5}\Omega(n)$, we have $\frac{1}{\sqrt{E_j}}<C_1 n^{-1}$, then
$$\sup_{x\in \mathbb{R}}|F_j(x)-\Phi_j(x)|<C_2\cdot n^{-1},$$
where $C_1$ and $C_2$ are constants which depend on $p$.

We claim an easy proposition, which is useful to generalize results in Berry-Esseen theorem for i.i.d. random vectors.
\begin{proposition}\label{control of independent sums}
Suppose $X_1,X_2,\ldots,X_n$, and $Y_1,Y_2,\ldots,Y_n$ are independent $d$-dimensional random vectors, $d\geq1$, with
cumulative distribution functions $F_1,F_2,\ldots,F_n$, and $G_1,G_2,\ldots,G_n$.
Let $$X=X_1+\ldots+X_n,$$
$$Y=Y_1+\ldots+Y_n,$$
and $F,G$ the cumulative distribution functions of $X,Y$.
If $$\sup_{x\in \mathbb{R}^d}|F_k(x)-G_k(x)|\leq\epsilon_k,$$
where $1\leq k\leq n$, then
$$\sup_{x\in \mathbb{R}^d}|F(x)-G(x)|\leq\sum_{i=1}^{n}\epsilon_i.$$
\end{proposition}

The proof of Proposition \ref{control of independent sums} can be given by induction on the number $n$ of summands in $X$ and $Y$.

\begin{comment}

\noindent{\sc Proof:}
It is sufficient to prove the proposition in the case $n=2$, because
the general case can be proved by induction on $n$.
For $n=2$, we have $F(x)=(F_1\ast F_2)(x)$£¬$G(x)=(G_1\ast G_2)(x)$,
where $A*B$ denote the convolution of functions $A$ and $B$.
Then for any $x\in \mathbb{R}_d$, we have
        \begin{align*}
            &|F(x)-G(x)|=|(F_1\ast F_2)(x)-(G_1\ast G_2)(x)|\\
            \leq&|(F_1\ast F_2)(x)-(G_1\ast F_2)(x)|+|(G_1\ast F_2)(x)-(G_1\ast G_2)(x)|\\
            =&\bigg|\int_{\mathbb{R}_d}(F_1(x-y)-G_1(x-y))dF_2(y)\bigg|+\bigg|\int_{\mathbb{R}_d}(F_2(x-y)-G_2(x-y))dG_1(y)\bigg|\\
            \leq&\int_{\mathbb{R}_d}|F_1(x-y)-G_1(x-y)|dF_2(y)+\int_{\mathbb{R}_d}|F_2(x-y)-G_2(x-y)|dG_1(y)\\
            \leq&\int_{\mathbb{R}_d}\epsilon_1 dF_2(y)+\int_{\mathbb{R}_d}\epsilon_2 dG_1(y)\\
            =&\epsilon_1+\epsilon_2.
        \end{align*}
$\Box$

\end{comment}

We go back to the proof of Lemma \ref{EX_1}. By Proposition \ref{control of independent sums},
there exist a constant $C_3$ which only depends on $p$, so that
$$\bigg|P(A_1;|i-pn|<n^{0.5}\Omega(n))-\Phi\bigg(\frac{i(n-2)p+(n-2i)}{\sigma}\bigg)\bigg|<C_3\cdot n^{-1}.$$
Since $\Phi(x)$ is continuous,
$$\lim_{n\rightarrow\infty} \Phi\bigg(\frac{i(n-2)p+(n-2i)}{\sigma}\bigg)=\Phi(1).$$
In addition, since $P(|i-pn|\geq n^{0.5}\Omega(n))<\epsilon_0$,
$$|P(A_1)-P(A_1;|i-pn|<n^{0.5}\Omega(n))|<\epsilon_0.$$
Since the choice of $\epsilon_0$ is arbitrary, we have
$$\lim_{n\rightarrow\infty}P(A_1)-\Phi(1)=0.$$
Lemma \ref{EX_1} is proved. $\Box$ \\

We will prove another lemma.
    \begin{lemma}\label{CovX_1,X_2}
        The limit $\lim_{n\rightarrow\infty} \textrm{Cov}[X_1,X_2]$ is equal to 0.
    \end{lemma}
We have changed the definition of the vertex sets $V_1,V_2$ in this part:
Let $V_1$ be the set of vertices which are adjacent to $v_1$, but not adjacent to $v_2$, and its size $|V_1|=i_1$;
$V_2$ the set of vertices which are adjacent to both $v_1$ and $v_2$, and its size $|V_2|=i_2$;
$V_3$ the set of vertices which are adjacent to $v_2$, but not adjacent to $v_1$, and its size $|V_3|=i_3$;
$V_4$ the set of vertices which are not adjacent to either $v_1$ or $v_2$, and its size $|V_4|=i_4=n-2-i_1-i_2-i_3$.
Then
    \begin{align*}
        &E[X_1 X_2]\\
        =&P(A_1\wedge A_2\wedge v_1\thicksim v_2)+P(A_1\wedge A_2\wedge v_1\nsim v_2)\\
        =&p\cdot P(A_1\wedge A_2|v_1\thicksim v_2)+(1-p)\cdot P(A_1\wedge A_2\wedge v_1\nsim v_2).
    \end{align*}
Denote $e_{jk}$,$1\leq j\leq k\leq 4$, be the number of edges whose two vertices are in $V_j$ and $V_k$, respectively.

\begin{comment}

Then
    \begin{displaymath}
        A_1\wedge A_2|v_1\sim v_2 \Leftrightarrow \left\{ \begin{aligned}
            &A_1|v_1\sim v_2\\
            &A_2|v_1\sim v_2
        \end{aligned} \right.
    \end{displaymath}
    \begin{displaymath}
        \Leftrightarrow \left\{ \begin{aligned}
            &2(i_1+i_2+1)\bigg(1+i_1+2i_2+i_3+\sum_{i\leq j\leq k\leq 4}e_{jk}\bigg)\\
            \leq &n(1+i_1+3i_2+i_3+2(e_{11}+e_{22}+e_{12})+(e_{13}+e_{14}+e_{23}+e_{24}))\\
            &\\
            &2(i_3+i_2+1)\bigg(1+i_1+2i_2+i_3+\sum_{i\leq j\leq k\leq 4}e_{jk}\bigg)\\
            \leq &n(1+i_1+3i_2+i_3+2(e_{33}+e_{22}+e_{23})+(e_{13}+e_{34}+e_{12}+e_{24})).
        \end{aligned} \right.
    \end{displaymath}
    \begin{displaymath}
        A_1\wedge A_2|v_1\nsim v_2 \Leftrightarrow \left\{ \begin{aligned}
            &A_1|v_1\nsim v_2\\
            &A_2|v_1\nsim v_2
        \end{aligned} \right.
    \end{displaymath}
    \begin{displaymath}
        \Leftrightarrow \left\{ \begin{aligned}
            &2(i_1+i_2)\bigg(i_1+2i_2+i_3+\sum_{i\leq j\leq k\leq 4}e_{jk}\bigg)\\
            \leq &n(i_1+2i_2+2(e_{11}+e_{22}+e_{12})+(e_{13}+e_{14}+e_{23}+e_{24}))\\
            &\\
            &2(i_3+i_2)\bigg(i_1+2i_2+i_3+\sum_{i\leq j\leq k\leq 4}e_{jk}\bigg)\\
            \leq &n(2i_2+i_3+2(e_{33}+e_{22}+e_{23})+(e_{13}+e_{34}+e_{12}+e_{24})).
        \end{aligned} \right.
    \end{displaymath}

\end{comment}

If $v_1$ is adjacent to $v_2$, the event $A_1\wedge A_2$ is equivalent to the following inequalities:
    \begin{displaymath}
        \left\{ \begin{aligned}
            B_1&\triangleq2(n-(i_1+i_2+1))(e_{11}+e_{22}+e_{12}+i_2)\\
            +&(n-2(i_1+i_2+1))(e_{13}+e_{14}+e_{23}+e_{24}+(i_1+i_2+i_3+1))\\
            -&2(i_1+i_2+1)(e_{33}+e_{34}+e_{44})\geq 0\\
            &\\
            B_2&\triangleq2(n-(i_3+i_2+1))(e_{33}+e_{22}+e_{23}+i_2)\\
            +&(n-2(i_3+i_2+1))(e_{13}+e_{34}+e_{12}+e_{24}+(i_1+i_2+i_3+1))\\
            -&2(i_3+i_2+1)(e_{11}+e_{14}+e_{44})\geq 0.
        \end{aligned} \right.
    \end{displaymath}

Since the $e_{jk}$'s are independent when the value of $i_1,i_2,i_3$ is fixed,
then the conditional expectation and variance of $B_1$ and $B_2$ are given by
    \begin{displaymath}
        E[B_1;i_1,i_2,i_3]=(n-(i_1+i_2+1))((n-2)p+(n-2(i_1+i_2+1))),
    \end{displaymath}
    \begin{displaymath}
        E[B_2;i_1,i_2,i_3]=(n-(i_3+i_2+1))((n-2)p+(n-2(i_3+i_2+1))),
    \end{displaymath}

    \begin{displaymath}
        \begin{aligned}
            &\textrm{Var}[B_1;i_1,i_2,i_3]\\
            =&\bigg[(2(n-(i_1+i_2+1)))^2\frac{(i_1+i_2+1)(i_1+i_2)}{2}\\
            &+(n-2(i_1+i_2+1))^2(i_1+i_2+1)(n-1-(i_1+i_2+1))\\
            &+(2(i_1+i_2+1))^2\frac{(n-1-(i_1+i_2+1))(n-2-(i_1+i_2+1))}{2}\bigg](p-p^2),\\
        \end{aligned}
    \end{displaymath}

    \begin{displaymath}
        \begin{aligned}
            &\textrm{Var}[B_2;i_1,i_2,i_3]\\
            =&\bigg[(2(n-(i_3+i_2+1)))^2\frac{(i_3+i_2+1)(i_3+i_2)}{2}\\
            &+(n-2(i_3+i_2+1))^2(i_3+i_2+1)(n-1-(i_3+i_2+1))\\
            &+(2(i_2+i_3+1))^2\frac{(n-1-(i_3+i_2+1))(n-2-(i_3+i_2+1))}{2}\bigg](p-p^2).
        \end{aligned}
    \end{displaymath}
In addition, the conditional covariance of $B_1$ and $B_2$ is given by
    \begin{displaymath}
        \begin{aligned}
            &\textrm{Cov}[B_1,B_2;i_1,i_2,i_3]\\
            =&\bigg[-4(n-(i_1+i_2+1))(i_2+i_3+1)\binom{i_1}{2}\\
            &+4(n-(i_1+i_2+1))(n-(i_2+i_3+1))\binom{i_2}{2}\\
            &-4(n-(i_2+i_3+1))(i_1+i_2+1)\binom{i_3}{2}\\
            &+4(i_1+i_2+1)(i_2+i_3+1)\binom{n-2-i_1-i_2-i_3}{2}\\
            &+2(n-(i_1+i_2+1))(n-2(i_2+i_3+1))i_1i_2\\
            &+(n-2(i_1+i_2+1))(n-2(i_2+i_3+1))i_1i_3\\
            &-2(n-2(i_1+i_2+1))(i_2+i_3+1)i_1(n-2-i_1-i_2-i_3)\\
            &+2(n-(i_2+i_3+1))(n-2(i_1+i_2+1))i_2i_3\\
            &+(n-2(i_1+i_2+1))(n-2(i_2+i_3+1))i_2(n-2-i_1-i_2-i_3)\\
            &-2(i_1+i_2+1)(n-2(i_2+i_3+1))i_3(n-2-i_1-i_2-i_3)\bigg](p-p^2).
        \end{aligned}
    \end{displaymath}

%% Similar to the case of $1$ dimension, since $i_1\thicksim B(n-2,p(1-p))$,
%% $i_2\thicksim B(n-2,p^2)$, $i_3\thicksim B(n-2,p(1-p))$,
%% then for any $\epsilon_1>0$, there exists $N\in \mathbb{N}$, so that for any $n>N$, we have
%% $$P(|i_1-np(1-p)|>n^{0.5}\Omega(n))<\epsilon_1,$$
%% $$P(|i_2-np^2|>n^{0.5}\Omega(n))<\epsilon_1,$$
%% and $$P(|i_3-np(1-p)|>n^{0.5}\Omega(n))<\epsilon_1.$$
Treat $B=(B_1,B_2)$ as a $2$-dimensional random vector.
Denote by $F(x_1,x_2)$ the cumulative distribution function of $\frac{B}{n^2}$,
$N=(N_1,N_2)$ the normal distribution with the same expectation and covariance matrix of $\frac{B}{n^2}$,
and $\Phi(x_1,x_2)$ the cumulative distribution function of $N$.

Note that the contribution of $e_{jk}$, $1\leq j\leq k\leq 4$, to $B_1,B_2$ are linear, by
Berry-Esseen theorem of independent identically distributed sequence on $\mathbb{R}^d$ proved by Bergstr{\"o}m \cite{B1945}\cite{C1948}, and Proposition \ref{control of independent sums},
If we have
    \begin{condition}\label{condition}
        $|i_1-np(1-p)|\leq n^{0.5}\Omega(n)$, $|i_2-np^2|\leq n^{0.5}\Omega(n)$, $|i_3-np(1-p)|\leq n^{0.5}\Omega(n)$,
    \end{condition}
    \noindent
then there exists a constant $C_3$ which only depends on $p$, so that
    \begin{displaymath}
        \sup|F(x_1,x_2)-\Phi(x_1,x_2)|<C_3\cdot n^{-1}.
    \end{displaymath}
When $n\rightarrow\infty$, and Condition \ref{condition} holds,
we have $$\textrm{Var}[B_j;i_1,i_2,i_3]=n^4 \textrm{Var}[N_j;i_1,i_2,i_3]=O(n^4), \quad j=1,2.$$
In addition, the covariance $$\textrm{Cov}[B_1,B_2;i_1,i_2,i_3]=n^4 \textrm{Cov}[N_1,N_2;i_1,i_2,i_3]=O(n^{3.5}\Omega(n))=o(n^4).$$
For any bounded region $R$ in $\mathbb{R}^2$
    \begin{displaymath}
        \lim_{n\rightarrow\infty,\delta_1\rightarrow 0}\sup_{(x_1,x_2)\in R}|P(N_1\leq x_1)P(N_2\leq x_2)-P(N_1\leq x_1,N_2\leq x_2)|=0.
    \end{displaymath}
Under Condition \ref{condition}, $\textrm{Var}[N_j;i_1,i_2,i_3]$ ($j=1,2$) are bounded, in addition,
when $n\rightarrow\infty$, the probability that Condition \ref{condition} holds approaches $1$. Thus,
    $$\lim_{n\rightarrow\infty}P(A_1,A_2|v_1\thicksim v_2)-P(A_1|v_1\thicksim v_2)P(A_2|v_1\thicksim v_2)=0.$$\\
When $v_1$ is not adjacent to $v_2$, by similar argument, we have
    $$\lim_{n\rightarrow\infty}P(A_1,A_2|v_1\nsim v_2)-P(A_1|v_1\nsim v_2)P(A_2|v_1\nsim v_2)=0.$$

Thus the equality $\lim_{n\rightarrow\infty} \textrm{Cov}[X_1,X_2]=0$ holds. Lemma \ref{CovX_1,X_2} is proved.\\

We go back to the final proof of Theorem \ref{main}, for any $\epsilon>0$,
    \begin{align*}
        &P\bigg(\bigg|\frac{N(n)}{n}-\Phi(1)\bigg|>\epsilon\bigg)\\
        \leq&P\bigg(\bigg|\frac{N(n)}{n}-E[X_1]\bigg|>\frac{\epsilon}{2}\bigg)+P\bigg(|E[X_1]-\Phi(1)|>\frac{\epsilon}{2}\bigg).
    \end{align*}
By Lemma \ref{EX_1}, when $n\rightarrow\infty$, the second term in the right side of the inequality above approaches $0$.
For the first term, by Chebyshev's Inequality,
    \begin{align*}
        &P\bigg(\bigg|\frac{N(n)}{n}-E[X_1]\bigg|>\frac{\epsilon}{2}\bigg)\\
        \leq&\frac{4\textrm{Var}[X_1+\ldots+X_n]}{n^2 \epsilon^2}.
    \end{align*}
In addition,
    \begin{align*}
        &\textrm{Var}[X_1+\ldots+X_n]\\
        =&\sum_{i=1}^{n}\textrm{Var}[X_1]+2\sum_{i<j}\textrm{Cov}[X_i,X_j]\\
        <&n+n(n-1)\textrm{Cov}[X_1,X_2]=o(n^2),
    \end{align*}
so the first term also approaches $0$ as $n\rightarrow\infty$.
We complete the proof of Theorem \ref{main}.

\section{In The Complete Bipartite Random Graph $B(n_1,n_2,p)$}\label{B(n_1,n_2,p)}

In the complete bipartite random graph $G=B(n_1,n_2,p)$, where $0<p<1$ is fixed, when $n_1$ and $n_2$ approaches infinity,
the asymptotic behavior of the proportion of groupies in $G$ depends on the ratio $\frac{n_1}{n_2}$.
Suppose this ratio $\frac{n_1}{n_2}$ has a limit $\alpha\geq0$ as $n_1,n_2\rightarrow\infty$.
%% i.e., $\lim_{n_1,n_2\rightarrow\infty}\frac{n_1}{n_2}=\alpha$.
Without loss of generality, we may assume that $\alpha\leq1$.
If $|n_1-n_2|\rightarrow\infty$, the proportion of groupies in $G$ converges in probability towards $\frac{1}{1+\alpha}$.

\begin{theorem}\label{Bipartite(ratio<1)}
    Suppose $G=B(n_1,n_2,p)$, $0<p<1$, is a complete bipartite random graph with $n=n_1+n_2$ vertices.
    Let $N(n_1,n_2)$ be the number of groupies in $G$. Then for any $\epsilon>0$, we have
    \begin{displaymath}
        P\bigg(\bigg|\frac{N(n_1,n_2)}{n}-\frac{1}{1+\alpha}\bigg|>\epsilon\bigg)\rightarrow0,
    \end{displaymath}
    as $n_1,n_2\rightarrow\infty$, $\frac{n_1}{n_2}\rightarrow\alpha\leq1$, and $|n_1-n_2|=c(n)\rightarrow\infty$.
\end{theorem}

Denote the two parts $P_1,P_2$ of vertices in $G$ by $v_{1,1},v_{1,2},\ldots,v_{1,n_1}$, and $v_{2,1},v_{2,2},\ldots,v_{2,n_2}$.
Let $A_{j,k}$ be the event that the vertex $v_{j,k}$ is a groupie, where $j=1,2$, and $1\leq k\leq n_{j}$;
$E_{j,k}$ the characteristic function of the event $A_{j,k}$, i.e., $E_{j,k}=1_{A_{j,k}}$.

\begin{comment}

Similar to the case of complete random graph $B(n,p)$, we may also exclude the case of the existence of an isolated vertex,
because

\begin{align*}
    & P(\textrm{there exists an isolated point in }B(n_1,n_2,p))\\
    \leq & n_1(1-p)^{n_2}+n_2(1-p)^{n_1}\\
    \rightarrow & 0, \\
    \textrm{ as }& n_1, n_2\rightarrow\infty \textrm{ and } \frac{n_1}{n_2}\rightarrow \alpha>0.
\end{align*}

\end{comment}

We claim that for vertices $v_{1,k}$ in the part $P_1$, the probability that $v_{1,k}$ is a groupie approaches 0,
and for vertices $v_{2,k}$ in the part $P_2$, the probability that $v_{2,k}$ is a groupie approaches 1.

\begin{lemma}\label{EX_1,1,EX_2,1(ratio<1)}
    As $n_1,n_2\rightarrow\infty$, if $n_1-n_2\rightarrow -\infty$, then
    \begin{displaymath}
            \lim E[X_{1,1}]=0,\\
    \end{displaymath}
    \begin{displaymath}
            \lim E[X_{2,1}]=1.
    \end{displaymath}
\end{lemma}

Note that $E[X_{j,1}]=P(A_{j,1})=P(A_{j,k})$, for any $j=1,2$, and $k=1,2,\ldots,n_1$.
We may first exclude the case that vertices $v_{1,1}$ or $v_{2,1}$ is isolated
because its probability $(1-p)^{n_2}$ or $(1-p)^{n_1}$ approaches $0$ as $n_1,n_2\rightarrow\infty$.
We then estimate the value of $E[X_{1,1}]$.
Suppose the vertex $v_{1,1}$ has degree $i$.
Let $V_1$ be the set of $i$ vertices that are adjacent to $v_{1,1}$,
and $V_2$ be the other $n_2-i$ vertices in the set $P_2$ that are not adjacent to $v_{1,1}$.
Let $e_j$ be the number of edges with one endpoint in the set $V_j$, where $j=1,2$,
and another endpoint other than $v_{1,1}$.
Then the event that $v_{1,1}$ is a groupie is equivalent to
\begin{displaymath}
    \frac{e_1+i}{i}\geq\frac{2(e_1+e_2+i)}{n_1+n_2},
\end{displaymath}
i.e.,
\begin{equation*}
    (n_1+n_2-2i)e_1-2ie_2+(n_1+n_2-2i)i\geq0.
\end{equation*}

Note that $e_1$ is subject to the binomial distribution $B(i(n_1-1),p)$,
and $e_2$ is subject to the binomial distribution $B((n_2-i)(n_1-1),p)$.
Then we have
\begin{align*}
    & E[(n_1+n_2-2i)e_1-2ie_2+(n_1+n_2-2i)i;i] \\
    = & (n_1+n_2-2i)i(n_1-1)p-2i(n_2-i)(n_1-1)p+(n_1+n_2-2i)i \\
    = & i((n_1-1)p(n_1-n_2)+(n_1+n_2-2i));
\end{align*}
and
\begin{align*}
    & \textrm{Var}[(n_1+n_2-2i)e_1-2ie_2+(n_1+n_2-2i)i;i]\\
    = & (n_1+n_2-2i)^2 i(n_1-1)p(1-p)-(2i)^2(n_2-i)(n_1-1)p(1-p) \\
    = & i(n_1-1)p(1-p)(n_1^2 + n_2^2 + 2 n_1 n_2 - 4 i n_1).
\end{align*}

Note that the degree $i$ of $v_{1,1}$ is subject to the binomial distribution $B(n_2,p)$.
Thus for any $\epsilon_0>0$, when $n_2$ approaches infinity,
we have $P(|i-n_2 p|\geq n_2^{0.5}\Omega(n_2))<\epsilon_0$.

Under the condition that $|i-n_2 p|<n_2^{0.5}\Omega(n_2)$, we have
\begin{displaymath}
    \mu\triangleq E[(n_1+n_2-2i)e_1-2ie_2+(n_1+n_2-2i)i;i]=O(n^2 c(n)),
\end{displaymath}
\begin{displaymath}
    \sigma^2\triangleq \textrm{Var}[(n_1+n_2-2i)e_1-2ie_2+(n_1+n_2-2i)i;i]=O(n^4).
\end{displaymath}
Note that $(n_1+n_2-2i)e_1-2ie_2+(n_1+n_2-2i)i$ is subject to the linear combination of binomial random variables.
In addition, as $n_1,n_2\rightarrow\infty$,
by the Berry-Esseen Theorem, there exists a constant $C$ which only depends on $p$, such that
\begin{displaymath}
    \bigg|P(A_{1,1})-\Phi\bigg(\frac{\mu}{\sigma}\bigg)\bigg|<C n^{-1}.
\end{displaymath}
As $n_1,n_2\rightarrow\infty$, we have $\frac{\mu}{\sigma}\rightarrow -\infty$, then
\begin{displaymath}
    \Phi\bigg(\frac{\mu}{\sigma}\bigg)\rightarrow 0.
\end{displaymath}
In addition, by $P(|i-p n_2|\geq n_2^{0.5}\Omega(n_2))<\epsilon_0$,
\begin{displaymath}
    |P(A_{1,1})-P(A_{1,1};P(|i-p n_2|\geq n_2^{0.5}\Omega(n_2))<\epsilon_0)|<\epsilon_0.
\end{displaymath}
For the arbitrary selection of $\epsilon_0>0$,

\begin{equation}\label{P(A_1,1)}
    \lim_{n\rightarrow\infty}P(A_{1,1})=0.
\end{equation}
By similar argument, we have

\begin{equation}\label{P(A_2,1)}
    \lim_{n\rightarrow\infty}P(A_{2,1})=1.
\end{equation}

The proof of Lemma \ref{EX_1,1,EX_2,1(ratio<1)} is done. $\Box$

\begin{remark}
    In Lemma \ref{EX_1,1,EX_2,1(ratio<1)} we do not need the condition that there is a limit $\alpha$ of $\frac{n_1}{n_2}$,
    as long as we have $n_1-n_2\rightarrow -\infty$.
\end{remark}

Then by (\ref{P(A_1,1)}) and (\ref{P(A_2,1)}),
\begin{displaymath}
    \lim_{n_1,n_2\rightarrow\infty,\frac{n_1}{n_2}\rightarrow\alpha,|n_1-n_2|\rightarrow\infty}\textrm{Cov}[X_{j_1,k_1},X_{j_2,k_2}]=0, \textrm{ where }(j_1,k_1)\neq(j_2,k_2).
\end{displaymath}

For any $\epsilon>0$,
    \begin{align*}
        &P\bigg(\bigg|\frac{N(n_1,n_2)}{n}-\frac{1}{1+\alpha}\bigg|>\epsilon\bigg)\\
        \leq&P\bigg(\bigg|\frac{N(n_1,n_2)}{n}-\frac{n_1}{n}E[X_{1,1}]-\frac{n_2}{n}E[X_{2,1}]\bigg|>\frac{\epsilon}{3}\bigg)\\
        &+P\bigg(\bigg|\frac{n_1}{n}E[X_{1,1}]\bigg|>\frac{\epsilon}{3}\bigg)+P\bigg(\bigg|\frac{n_2}{n}E[X_{2,1}]-\frac{1}{1+\alpha}\bigg|>\frac{\epsilon}{3}\bigg).
    \end{align*}
By Lemma \ref{EX_1,1,EX_2,1(ratio<1)}, when $n_1,n_2\rightarrow\infty$, and $\frac{n_1}{n_2}\rightarrow\alpha$,
the second and third terms in the right side of the inequality above approaches to $0$.
For the first term, by Chebyshev's Inequality,
    \begin{align*}
        &P\bigg(\bigg|\frac{N(n_1,n_2)}{n}-\frac{n_1}{n}E[X_{1,1}]-\frac{n_2}{n}E[X_{2,1}]\bigg|>\frac{\epsilon}{3}\bigg)\\
        \leq&\frac{9\textrm{Var}[X_{1,1}+\ldots+X_{1,n_1}+X_{2,1}+\ldots+X_{2,n_2}]}{n^2 \epsilon^2}.
    \end{align*}
In addition,
    \begin{align*}
        &\textrm{Var}[X_{1,1}+\ldots+X_{1,n_1}+X_{2,1}+\ldots+X_{2,n_2}]\\
        =&\sum_{k=1}^{n_1}\textrm{Var}[X_{1,k}]+\sum_{k=1}^{n_2}\textrm{Var}[X_{2,k}]+2\sum_{(j_1,k_1)\neq(j_2,k_2)}\textrm{Cov}[X_{j_1,k_1},X_{j_2,k_2}]\\
        <&n+n(n-1)\textrm{Cov}[X_1,X_2]=o(n^2),
    \end{align*}
so the first term also approaches to $0$ as $n\rightarrow\infty$.
This completes the proof of Theorem \ref{Bipartite(ratio<1)}. $\Box$ \\

For the case of $\alpha=1$, we claim that the limit of proportion of groupies can be determined by the difference $c$ of the number of vertices in the two parts.
By the symmetry of the two parts, we may assume $c=n_1-n_2>0$. In the following argument, let $n_1$ and $n_2$ approaches infinity with the difference $c=n_1-n_2$ fixed, then the proportion of groupies in bipartite graph $B(n_1,n_2,p)$ converges in probability towards $\frac{1}{2}(\Phi(1+\frac{pc}{2(1-p)})+\Phi(1-\frac{pc}{2(1-p)}))$.

\begin{theorem}\label{Bipartite(ratio=1)}
    Suppose $G=B(n_1,n_2,p)$, $0<p<1$, is a complete bipartite random graph with $n=n_1+n_2$ vertices.
    Let $N(n_1,n_2)$ be the number of groupies in $G$.
    Then for any $\epsilon>0$, and $c \in \mathbb{Z}$, we have
    \begin{align*}
        &P\bigg(\bigg|\frac{N(n_1,n_2)}{n}-\frac{1}{2}\bigg(\Phi\bigg(1+\frac{pc}{2(1-p)}\bigg)+\Phi\bigg(1-\frac{pc}{2(1-p)}\bigg)\bigg)\bigg|
        >\epsilon\bigg)\rightarrow 0 ,\\
        &\textrm{as }n_1,n_2\rightarrow\infty, \textrm{ and }n_1-n_2=c.
    \end{align*}
\end{theorem}

We use the same notations in the proof of Theorem \ref{Bipartite(ratio<1)}.
As the number of vertices goes to infinite, the probability that the vertex $v_{1,k}$ (or $v_{2,k}$) in part $P_1$ (or $P_2$) approaches to $\Phi(1+\frac{pc}{2(1-p)})$ (or $\Phi(1-\frac{pc}{2(1-p)})$).

\begin{lemma}\label{EX_1,1}
    \begin{displaymath}
        \lim_{n_1,n_2\rightarrow\infty,n_1-n_2=c}E[X_{1,1}]=\Phi\bigg(1+\frac{pc}{2(1-p)}\bigg),\\
    \end{displaymath}
    \begin{displaymath}
        \lim_{n_1,n_2\rightarrow\infty,n_1-n_2=c}E[X_{2,1}]=\Phi\bigg(1-\frac{pc}{2(1-p)}\bigg).
    \end{displaymath}
\end{lemma}

The proofs of the two equalities are similar, so we just give the proof of the first one.
The event that the vertex $v_{1,1}$ is a groupie, is equivalent to
\begin{displaymath}
    (n_1+n_2-2i)e_1-2i e_2+(n_1+n_2-2i)i \geq 0.
\end{displaymath}

Given the condition $n_1-n_2=c$ and the value of $i$, the conditional expectation
\begin{align*}
    \mu\triangleq & E[(n_1+n_2-2i)e_1-2ie_2+(n_1+n_2-2i)i;i] \\
    = & i((n_1-1)p(n_1-n_2)+(n_1+n_2-2i)) \\
    = & i((n_2+c-1)pc+(2n_2+c-2i));
\end{align*}
and
\begin{align*}
    \sigma^2\triangleq & \textrm{Var}[(n_1+n_2-2i)e_1-2ie_2+(n_1+n_2-2i)i;i] \\
    = & i(n_1-1)p(1-p)(n_1^2 + n_2^2 + 2 n_1 n_2 - 4 i n_1).
\end{align*}

Fix the value of $i$ under the condition $|i-n_2 p|<n_2^{0.5}\Omega(n_2)$,
by the 1-dimentional Berry-Esseen theorem and Proposition \ref{control of independent sums}, we have
\begin{displaymath}
    |P(v_{1,1}\textrm{ is a groupie};|i-n_2 p|<n_2^{0.5}\Omega(n_2))-\Phi(\frac{\mu}{\sigma})|=O(n^{-1}).
\end{displaymath}
Note that as $n_1,n_2\rightarrow\infty$, and $n_1-n_2=c$,
\begin{displaymath}
    \lim \Phi\bigg(\frac{\mu}{\sigma}\bigg)=\Phi\bigg(1+\frac{pc}{2(1-p)}\bigg).
\end{displaymath}
In addition, the condition $|i-n_2 p|<n_2^{0.5}\Omega(n_2)$ holds with probability $1-o(1)$. Thus
\begin{displaymath}
    \lim_{n_1,n_2\rightarrow\infty,n_1-n_2=c}E[X_{1,1}]=\Phi\bigg(1+\frac{pc}{2(1-p)}\bigg)
\end{displaymath}
holds. $\Box$

For the covariance of the random variables $X_{j,k}$, where $j=1,2$ and $1\leq k\leq n_j$,
it approaches $0$, as $n_1,n_2\rightarrow\infty$ and $n_1-n_2=c$.

\begin{lemma}\label{CovX_j_1,k_1,X_j_2,k_2}
    As $n_1,n_2\rightarrow\infty$ and $n_1-n_2=c$,
    the covariance of characteristic functions of two distinct events $A_{j,k}$ approaches 0,
    i.e., for any two distinct pairs $(j_1,k_1)$ and $(j_2,k_2)$,
    \begin{equation*}
        \lim \textrm{Cov}[X_{j_1,k_1},X_{j_2,k_2}]=0.
    \end{equation*}
\end{lemma}

By the symmetry of the complete bipartite random graph $B(n_1,n_2,p)$, we only need to show that
\begin{equation}\label{CovX_1,1,X_1,2}
    \lim \textrm{Cov}[X_{1,1},X_{1,2}]=0,
\end{equation}
and
\begin{equation}\label{CovX_1,1,X_2,1}
    \lim \textrm{Cov}[X_{1,1},X_{2,1}]=0.
\end{equation}

For the proof of equality (\ref{CovX_1,1,X_1,2}), the notations are defined as follows.
Let $V_1$ be the set of vertices in part $P_2$ which are adjacent to $v_{1,1}$, but not adjacent to $v_{1,2}$;
$V_2$ the set of vertices in part $P_2$ which are adjacent to both $v_{1,1}$ and $v_{1,2}$;
$V_3$ the set of vertices in part $P_2$ which are adjacent to $v_{1,2}$, but not adjacent to $v_{1,1}$;
$V_4$ the set of vertices in part $P_2$ which are not adjacent to either $v_{1,1}$ or $v_{1,2}$.
Let $i_j$ be the size $|V_{j}|$, and $e_j$ the number of edges between part $P_1$ (other than $v_{1,1},v_{1,2}$) and vertex set $V_{j}$ in $P_2$,
where $j=1,2,3,4$. Then there is a relation $i_1+i_2+i_3+i_4=n_2$.

Thus the event $A_{1,1}$ that the vertex $v_{1,1}$ is a groupie is equivalent to
\begin{displaymath}
    \frac{e_1+e_2+i_1+2i_2}{i_1+i_2}\geq\frac{2(i_1+2i_2+i_3+e_1+e_2+e_3+e_4)}{n_1+n_2}
\end{displaymath}
i.e.,
\begin{align*}
    B_1 \triangleq &(n_1+n_2-2(i_1+i_2))(e_1+e_2)-2(i_1+i_2)(e_3+e_4)\\
    &+(i_1+i_2)(n_1+n_2-2(i_1+2i_2+i_3))+i_2(n_1+n_2) \geq 0.
\end{align*}
Similarly, the event $A_{1,2}$ is equivalent to
\begin{align*}
    B_2 \triangleq &(n_1+n_2-2(i_2+i_3))(e_2+e_3)-2(i_2+i_3)(e_1+e_4)\\
    &+(i_2+i_3)(n_1+n_2-2(i_1+2i_2+i_3))+i_2(n_1+n_2) \geq 0.
\end{align*}

Note that when $i_1,i_2,i_3$ (and $i_4$) is fixed, $e_j$ is subject to binomial distribution $B(i_j(n_1-2),p)$, where $j=1,2,3,4$.
Applying the relation that $n_1-n_2=c$ and by brute force, the conditional expectations of $B_1$ and $B_2$ are given by
\begin{align*}
    E[B_1;i_1,i_2,i_3]=&(i_1+i_2)((n_1-2)pc+n_1+n_2-2(i_1+2i_2+i_3))+i_2(n_1+n_2)\\
    E[B_2;i_1,i_2,i_3]=&(i_2+i_3)((n_1-2)pc+n_1+n_2-2(i_1+2i_2+i_3))+i_2(n_1+n_2).
\end{align*}
The conditional variance and covariance of $B_1$ and $B_2$ are given by
\begin{align*}
    \textrm{Var}[B_1;i_1,i_2,i_3]=&(n_1+n_2-2(i_1+i_2))^2(i_1+i_2)(n_1-2)p(1-p)\\
    &+4(i_1+i_2)^2(i_3+i_4)(n_1-2)p(1-p)\\
    \textrm{Var}[B_2;i_1,i_2,i_3]=&(n_1+n_2-2(i_2+i_3))^2(i_2+i_3)(n_1-2)p(1-p)\\
    &+4(i_2+i_3)^2(i_1+i_4)(n_1-2)p(1-p)\\
    \textrm{Cov}[B_1,B_2;i_1,i_2,i_3]=&((n_1+n_2-2(i_1+i_2))(n_1+n_2-2(i_2+i_3))i_2\\
    &-2(n_1+n_2-2(i_1+i_2))(i_2+i_3)i_1\\
    &-2(n_1+n_2-2(i_2+i_3))(i_1+i_2)i_3\\
    &+4(i_1+i_2)(i_2+i_3)i_4)(n_1-2)p(1-p).
\end{align*}

Analogue to the argument for the complete random graph $B(n,p)$, note that $i_1\sim B(n_2,p(1-p))$,
$i_2\sim B(n_2,p^2)$, and $i_3\sim B(n_2,p(1-p))$. Then the condition
\begin{condition}
    $|i_1-n_2 p(1-p)|\leq n_2^{0.5}\Omega(n_2)$, $|i_2-n_2 p^2|\leq n_2^{0.5}\Omega(n_2)$, $|i_3-n_2 p(1-p)|\leq n_2^{0.5}\Omega(n_2)$
\end{condition}
holds with probability $1-o(1)$ as $n_1,n_2\rightarrow\infty$, and $n_1-n_2=c$.
Under this condition, $$E[B_1;i_1,i_2,i_3]=O(n^2),$$
$$E[B_2;i_1,i_2,i_3]=O(n^2),$$
$$\textrm{Var}[B_1;i_1,i_2,i_3]=O(n^4),$$
$$\textrm{Var}[B_2;i_1,i_2,i_3]=O(n^4),$$
and $$\textrm{Cov}[B_1,B_2;i_1,i_2,i_3]=o(n^4).$$
By the same argument for complete random graph $B(n,p)$,
we use the normal distribution with the same expectation and covariance matrix of $(\frac{B_1}{n^2},\frac{B_2}{n^2})$,
to estimate the random vector $(\frac{B_1}{n^2},\frac{B_2}{n^2})$, and finally get $$\lim_{n_1,n_2\rightarrow\infty,n_1-n_2=c}\textrm{Cov}[X_{1,1},X_{1,2}]=0.$$

For the proof of equality (\ref{CovX_1,1,X_2,1}), the notations are defined as follows:
Let $V_1$ be the set of vertices other than $v_{2,1}$ in part $P_2$ which are adjacent to $v_{1,1}$;
$V_2$ the set of vertices other than $v_{1,1}$ in part $P_1$ which are adjacent to $v_{2,1}$;
$V_3$ the set of vertices other than $v_{2,1}$ in part $P_2$ which are not adjacent to $v_{1,1}$;
$V_4$ the set of vertices other than $v_{1,1}$ in part $P_1$ which are not adjacent to $v_{2,1}$.
Let $i_j$ be the size $|V_j|$ and $e_{j,k}$ the number of edges between vertex set $V_j$ and $V_k$, where $j,k=1,2,3,4$.
Then there are relations $i_1+i_3=n_2-1$, and $i_2+i_4=n_1-1$.
If $v_{1,1}$ and $v_{2,1}$ are not joined by an edge, then the event ``$v_{1,1}$ is a groupie'' is equivalent to
\begin{align*}
    B_1 \triangleq & (n_1+n_2-i_1)(e_{1,2}+e_{1,4})-i_1(e_{3,2}+e_{3,4})\\
    & +i_1(n_1+n_2-i_1-i_2) \geq 0,\\
\end{align*}
and the event ``$v_{2,1}$ is a groupie'' is equivalent to
\begin{align*}
    B_2 \triangleq & (n_1+n_2-i_2)(e_{1,2}+e_{3,2})-i_2(e_{1,4}+e_{3,4})\\
    & +i_2(n_1+n_2-i_1-i_2) \geq 0.\\
\end{align*}

Analogue to the argument for the equality (\ref{CovX_1,1,X_1,2}), note that $i_1\sim B(n_2-1,p)$, and
$i_2\sim B(n_1-1,p)$. Then the condition
\begin{condition}
    $|i_1-(n_2-1)p|\leq n_2^{0.5}\Omega(n_2)$, $|i_2-(n_1-1)p|\leq n_1^{0.5}\Omega(n_1)$
\end{condition}
holds with probability $1-o(1)$ as $\delta_1\rightarrow 0$, $n_1,n_2\rightarrow\infty$, and $n_1-n_2=c$.
Under this condition, by brute force, we have $$E[B_1;i_1,i_2]=O(n^2),$$
$$E[B_2;i_1,i_2]=O(n^2),$$
$$\textrm{Var}[B_1;i_1,i_2]=O(n^4),$$
$$\textrm{Var}[B_2;i_1,i_2]=O(n^4),$$
and $$\textrm{Cov}[B_1,B_2;i_1,i_2]=o(n^4).$$
By the same argument for complete random graph $B(n,p)$,
we use the normal distribution with the same expectation and covariance matrix of $(\frac{B_1}{n^2},\frac{B_2}{n^2})$,
to estimate the random vector $(\frac{B_1}{n^2},\frac{B_2}{n^2})$, and finally get $$\lim_{n_1,n_2\rightarrow\infty,n_1-n_2=c}\textrm{Cov}[X_{1,1},X_{2,1}]=0.$$ under the assumption that $v_{1,1}$ and $v_{2,1}$ are not adjacent.
If $v_{1,1}$ and $v_{2,1}$ are adjacent to each other, the argument is similar. So the equality (\ref{CovX_1,1,X_2,1}) holds.
$\Box$

The remaining argument is similar to the case of $\frac{n_1}{n_2}\rightarrow\alpha$ and $|n_1-n_2|\rightarrow\infty$.
For any $\epsilon>0$,
    \begin{align*}
        &P\bigg(\bigg|\frac{N(n_1,n_2)}{n}-\frac{1}{2}\bigg(\Phi\bigg(1+\frac{pc}{2(1-p)}\bigg)+\Phi\bigg(1-\frac{pc}{2(1-p)}\bigg)\bigg)\bigg|>\epsilon\bigg)\\
        \leq & P\bigg(\bigg|\frac{N(n_1,n_2)}{n}-\frac{n_1}{n}E[X_{1,1}]-\frac{n_2}{n}E[X_{2,1}]\bigg|>\frac{\epsilon}{3}\bigg)\\
        &+P\bigg(\bigg|\frac{n_1}{n}E[X_{1,1}]-\frac{1}{2}\bigg(\Phi\bigg(1+\frac{pc}{2(1-p)}\bigg)\bigg)\bigg|>\frac{\epsilon}{3}\bigg)\\
        &+P\bigg(\bigg|\frac{n_2}{n}E[X_{2,1}]-\frac{1}{2}\bigg(\Phi\bigg(1-\frac{pc}{2(1-p)}\bigg)\bigg)\bigg|>\frac{\epsilon}{3}\bigg).
    \end{align*}
By Lemma \ref{EX_1,1}, when $n_1,n_2\rightarrow\infty$, and $\frac{n_1}{n_2}\rightarrow\alpha$,
the second and third terms in the right side of the inequality above approaches to $0$.
For the first term, by Chebyshev's Inequality,
    \begin{align*}
        &P\bigg(\bigg|\frac{N(n_1,n_2)}{n}-\frac{1}{2}\bigg(\Phi\bigg(1+\frac{pc}{2(1-p)}\bigg)+\Phi\bigg(1-\frac{pc}{2(1-p)}\bigg)\bigg)\bigg|>\frac{\epsilon}{3}\bigg)\\
        \leq & \frac{9\textrm{Var}[X_{1,1}+\ldots+X_{1,n_1}+X_{2,1}+\ldots+X_{2,n_2}]}{n^2 \epsilon^2}.
    \end{align*}
In addition,
    \begin{align*}
        & \textrm{Var}[X_{1,1}+\ldots+X_{1,n_1}+X_{2,1}+\ldots+X_{2,n_2}]\\
        =&\sum_{k=1}^{n_1}\textrm{Var}[X_{1,k}]+\sum_{k=1}^{n_2}\textrm{Var}[X_{2,k}]+2\sum_{(j_1,k_1)\neq(j_2,k_2)}\textrm{Cov}[X_{j_1,k_1},X_{j_2,k_2}]\\
        <&n+n(n-1)(\textrm{Cov}[X_{1,1},X_{1,2}]+\textrm{Cov}[X_{1,1},X_{2,1}])=o(n^2),
    \end{align*}
so the first term also approaches to $0$ as $n\rightarrow\infty$.
This completes the proof of Theorem \ref{Bipartite(ratio=1)}. $\Box$

If there is no limit for the ratio $\frac{n_1}{n_2}$ of vertices in the two parts,
the proportion of groupies in $B(n_1,n_2,p)$ does not converge in probability, known from Lemma \ref{EX_1,1,EX_2,1(ratio<1)}.
By Theorem \ref{Bipartite(ratio<1)} and Theorem \ref{Bipartite(ratio=1)},
we determine whether the proportion of groupies in complete bipartite random graphs $B(n_1,n_2,p)$ converges in probability
as the number of vertices $n_1$ and $n_2$ in the two parts both approaches infinity in general,
and give the limit if it converges in probability.

\section*{Acknowledgements}
The author is grateful to Professor Chunwei Song for his useful comments and suggestions.

\bibliography{pkuthss}

\begin{thebibliography}{10}

\bibitem{AKS1980}
Mikl{\'o}s Ajtai, J{\'a}nos Koml{\'o}s, and Endre Szemer{\'e}di.
\newblock A note on {R}amsey numbers.
\newblock {\em J. Combin. Theory Ser. A}, 29(3):354--360, 1980.

\bibitem{B1945}
Harald Bergstr{\"o}m.
\newblock On the central limit theorem in the space {$R_k, k>1$}.
\newblock {\em Skand. Aktuarietidskr.}, 28:106--127, 1945.

\bibitem{B1941}
Andrew~C. Berry.
\newblock The accuracy of the {G}aussian approximation to the sum of
  independent variates.
\newblock {\em Trans. Amer. Math. Soc.}, 49:122--136, 1941.

\bibitem{BEHST1994}
E.~Bertram, P.~Erd{\H{o}}s, P.~Hor{\'a}k, J.~{\v{S}}ir{\'a}{\v{n}}, and Zs.
  Tuza.
\newblock Local and global average degree in graphs and multigraphs.
\newblock {\em J. Graph Theory}, 18(7):647--661, 1994.

\bibitem{C1948}
Kai~Lai Chung.
\newblock Asymptotic distribution of the maximum cumulative sum of independent
  random variables.
\newblock {\em Bull. Amer. Math. Soc.}, 54:1162--1170, 1948.

\bibitem{VT2009}
W.~Fernandez de~la Vega and Zsolt Tuza.
\newblock Groupies in random graphs.
\newblock {\em Inform. Process. Lett.}, 109(7):339--340, 2009.

\bibitem{E1942}
Carl-Gustav Esseen.
\newblock On the {L}iapounoff limit of error in the theory of probability.
\newblock {\em Ark. Mat. Astr. Fys.}, 28A(9):19, 1942.

\bibitem{H2007}
Pak~Tung Ho.
\newblock On groupies in graphs.
\newblock {\em Australas. J. Combin.}, 38:173--177, 2007.

\bibitem{M1996}
John Mackey.
\newblock A lower bound for groupies in graphs.
\newblock {\em J. Graph Theory}, 21(3):323--326, 1996.

\bibitem{PST1995}
Svatopluk Poljak, Tam{\'a}s Szab{\'o}, and Zsolt Tuza.
\newblock Extremum and convergence of local average degrees in graphs.
\newblock In {\em Proceedings of the {T}wenty-sixth {S}outheastern
  {I}nternational {C}onference on {C}ombinatorics, {G}raph {T}heory and
  {C}omputing ({B}oca {R}aton, {FL}, 1995)}, volume 112, pages 191--198, 1995.

\bibitem{S2010}
Yilun Shang.
\newblock Groupies in random bipartite graphs.
\newblock {\em Appl. Anal. Discrete Math.}, 4(2):278--283, 2010.

\end{thebibliography}

\end{document}